\documentclass[12pt]{article}

\usepackage{graphicx, epsfig}
\usepackage{amsmath,amssymb}
\usepackage[authoryear,nonamebreak,elide]{natbib}
\usepackage[latin1]{inputenc}

\usepackage{color}
\usepackage[pagebackref,colorlinks=true,urlcolor=red,citecolor=blue,linkcolor=red]{hyperref}
\usepackage{bibentry}

\begin{document}

\title{Small-sample one-sided testing in extreme value regression models}
\author{Silvia L.P. Ferrari\thanks{Corresponding author. Email: silviaferrari.usp@gmail.com}
\ \ and \ \  Eliane C. Pinheiro\\
{\small {\em Departamento de Estat\' \i stica, Universidade de S\~ao Paulo,
Brazil}}
}

\maketitle
\begin{abstract}

We derive adjusted signed likelihood ratio statistics for a general class of extreme value regression models.
The adjustments reduce the error in the standard normal approximation to the distribution of the signed likelihood ratio 
statistic. 
We use Monte Carlo simulations to compare the finite-sample performance of the different tests.
Our simulations suggest that the signed likelihood ratio test tends to be liberal when the sample size is not large,
and that the adjustments are effective in shrinking the size distortion. Two real data applications are presented
and discussed.


\end{abstract}

\section{Introduction}
\label{introduction}

Statistics of extremes are applied in a wide range of areas, such as marine, civil and telecommunications  engineering, metallurgy, geology, astronomy, oceanography, meteorology, hydrology, reliability, finance, insurance, medicine and nutrition. 
In the last few decades, the growth in the intensity and frequency of environmental phenomena -- Hurricane Katrina (2005), 
the Japan earthquake (2011), mudslides and floods due to intense rains in Australia (2010-2011), Pakistan (2010) and Brazil (2011), 
the typhoon in the Philippines (2013) -- and disasters caused by man -- oil leaks (2010) and financial crises (2008, 2010) -- has highlighted the necessity of improving the inference accuracy of the  statistics of extremes.

A classical reference regarding the statistics of extremes is the study \cite{GUMBEL1958}, who formalized the extreme value theory 
presenting the only three asymptotic limiting distributions of extreme values obtained by \cite{FRECHET} and \cite{FISHER}, 
and unified by \cite{GNEDENKO}. More recent book-length references include the works of 
\cite{HAAN}, 
\cite{COLES}, 
\cite{BEIRLANT}, 
\cite{KOTZ2000} 
and, 
\cite{CASTILLO}, 
among others. 

The extreme value distribution is frequently used to model extreme events data and reliability/survival outcomes. 
The focus of this paper is on inference for extreme value regression models when a sample is small or of
moderate size. Specifically, we deal with a general class of extreme value regression models introduced by \cite{BARRETO} in which location and dispersion parameters are linear or, possibly, non-linear functions of  regressors. 
\cite{MUELLERRUFIBACH} commented on the lack of small sample extreme value theory in the literature.
Nevertheless, there are many practical applications for which data sets are not large enough to safely apply asymptotic theory.
Recent contributions in small-sample inference for extreme value regression models are those made by \cite{BARRETO}, who obtained bias corrections for
maximum likelihood estimators, and \cite{FERRARI2012}, who derived Skovgaard's adjustment \citep{SKOVGAARD2001} to the
likelihood ratio statistic for testing two-sided hypotheses. \cite{FERRARI2012} compared the small sample accuracy of the usual likelihood ratio, Wald, score and gradient tests, and the adjusted likelihood ratio tests in extreme value regression models with linear and non-linear predictors for both location and dispersion parameters. Simulation results suggest that the type I error of the likelihood ratio and Wald tests can be markedly greater than the nominal level in small- and moderate-sized samples. The gradient test presents similar behavior, but is much less oversized. The score test is even less size distorted, can be conservative in some scenarios, and is competitive with the adjusted likelihood ratio test in some cases. Overall, the adjusted likelihood ratio test performs better than all others tests.

The present paper is focused on one-sided tests, for which the above mentioned statistics are not suitable. 
The signed likelihood ratio statistic is the most widely used for one-sided testing of a simple null hypothesis about a scalar parameter in the presence of nuisance parameters. This statistic has a limiting standard normal distribution, but its finite-sample law
is usually unknown. It is well-known that the use of the asymptotic distribution in small samples frequently lead to 
oversized signed likelihood ratio tests. Considerable research has been conducted to propose adjustments to improve the accuracy 
of the standard normal approximation to the distribution of the signed root of the likelihood ratio statistic. In this paper, we derive five adjustments, proposed by \citet{BARNDORFF-NIELSEN1986,BARNDORFF-NIELSEN1991},  \cite{DICICCIO}, \cite{SKOVGAARD1996}, \cite{SEVERINI1999} and \cite{FRASER1999}, to extreme value regression models. We compare the finite-sample performance of the signed likelihood ratio test and the adjusted signed likelihood ratio tests obtained in this work. Simulation results suggest that the adjusted statistics approximately follow a standard normal distribution with a high degree of accuracy even in small samples. Furthermore, we illustrate an application of the usual signed likelihood ratio test and its modified versions in a real data set.

The paper is organized as follows. In Section~\ref{model} we present the extreme value regression models.
In Section~\ref{aslrs}, the modified signed likelihood ratio statistic proposed by
\citet{BARNDORFF-NIELSEN1986,BARNDORFF-NIELSEN1991} is presented along with approximations derived by different authors, namely \cite{DICICCIO}, \cite{SKOVGAARD1996}, \cite{SEVERINI1999} and \cite{FRASER1999}. In Section~\ref{aslrsevrm} we derive adjusted signed likelihood ratio statistics in extreme value regression models. Monte Carlo simulations are presented in Section~\ref{MonteCarlo} to compare the performances of the tests, whereas applications to real data sets are presented in Section~\ref{applications}. 
The paper ends with conclusions in Section~\ref{conclusion}. Technical details regarding DiCiccio \& Martin's
approach are presented in the Appendix.

\section{Extreme value regression models}
\label{model}

Let $y$ be a continuous random variable with a maximum extreme value distribution with location parameter $\mu$, dispersion parameter $\sigma$, and density function 
\begin{equation}
\label{extreme value-density}
f(y; \mu, \sigma) = \frac{1}{\sigma}\exp\left(-\frac{y-\mu}{\sigma}\right)\exp\left\{-\exp\left(-\frac{y-\mu}{\sigma}\right)\right\},\quad y \in I\!\! R,
\end{equation}
where $\mu \in I\!\! R$ and $\sigma > 0$;  we write $y\sim EV_{max}(\mu,\sigma)$.
The mean and the variance of $y$ are
${\rm E}(y) = \mu+{\mathcal E}\sigma$
and 
${\rm var}(y) = {\sigma^2\pi^2}/{6},$
respectively, where ${\mathcal E}$ is the Euler constant; ${\mathcal E}\approx 0.5772$.
This distribution is also called the Gumbel or type I extreme value distribution. Here, we refer to distribution as the maximum 
extreme value distribution to distinguish it from the minimum extreme value distribution, 
which is also often called the Gumbel or type I extreme value distribution in the statistical literature. 
The distribution in (\ref{extreme value-density}) is one of the three possible limiting laws 
of the standardized maximum of independent and identically distributed random variables \citep{GNEDENKO},
and is frequently invoked to model extreme events; see, for instance, \citet[Tables 9.16]{CASTILLO} and \citet[Section 3.4.1]{COLES}.

The density function of the minimum extreme value distribution with location parameter $\mu$ and dispersion parameter $\sigma$ is given by
\begin{equation}
\label{eq:min-extreme value-density}
f(y; \mu, \sigma) = \frac{1}{\sigma}\exp\left(\frac{y-\mu}{\sigma}\right)\exp\left\{-\exp\left(\frac{y-\mu}{\sigma}\right)\right\} ,\quad y \in I\!\! R,
\end{equation}
and we write $y\sim EV_{min}(\mu,\sigma)$. The mean of the minimum extreme value distribution is $\mu-{\mathcal E}\sigma$
and its variance coincides with that of the maximum extreme value distribution.
A useful property is that
\begin{eqnarray}\label{min-max}
y\sim EV_{min}(\mu,\sigma) \Longleftrightarrow -y\sim EV_{max}(-\mu,\sigma).
\end{eqnarray}
The distribution in (\ref{eq:min-extreme value-density}) is the distribution of the logarithm of Weibull distributed random variables,
and is often used in reliability and survival analysis to model log-lifetimes \citep{LAWLESS}.

Let $y_1, \ldots, y_n$ be independent random variables, where 
\begin{equation}
\label{extreme value-max}
y_t\sim EV_{max}(\mu_t,\sigma_t), 
\end{equation}
for $t=1,\ldots,n$. The maximum extreme value regression model is defined by (\ref{extreme value-max}) and by two systematic components given by
\begin{equation}\label{linkmu}
g(\mu_t) = \eta_t = \eta(x_{t},\beta)
\end{equation}
and
\begin{equation}\label{linkphi}
h(\sigma_t) =  \delta_t = \delta(z_{t},\gamma),
\end{equation}
where $\beta = (\beta_1, \ldots, \beta_k)^{\!\top}$ and $\gamma = (\gamma_1, \ldots, \gamma_m)^{\!\top}$ are vectors of  unknown regression parameters ($\beta \in I\!\! R^k$ and $\gamma \in I\!\! R^m$, $k + m < n$) and $x_{t}$ and $z_{t}$ are observations on covariates. Here, $\eta(\cdot,\cdot)$ and $\delta(\cdot,\cdot)$ are continuously twice differentiable 
(possibly nonlinear) functions in the second argument.
Finally, $g(\cdot)$ and $h(\cdot)$ are known, strictly monotonic and twice differentiable link functions that map $I\!\! R$ and $I\!\! R^+$ 
onto $I\!\! R$, respectively. Let $X$ be the derivative matrix of $\eta=(\eta_1,\ldots,\eta_n)^{\top}$ with respect to $\beta^\top$, and let $Z$ be the derivative matrix of $\delta=(\delta_1,\ldots,\delta_n)^{\top}$ with respect to $\gamma^\top$. It is assumed that ${\rm rank}(X)=k$ and ${\rm rank}(Z)=m$ for all $\beta$ and $\gamma$.

Analogously, we define the minimum extreme value regression model as above, by substituting the following for (\ref{extreme value-max}):
\begin{equation}
\label{extreme value-min}
y_t\sim EV_{min}(\mu_t,\sigma_t), 
\end{equation}
for $t=1,\ldots,n$. Hereinafter, we shall refer to the maximum extreme value regression model (\ref{extreme value-max})--(\ref{linkphi}) only,
but, as will be shown, our results can be easily adapted to the minimum extreme value regression model  
(\ref{linkmu})--(\ref{extreme value-min}).

\section{Adjusted signed likelihood ratio statistics}\label{aslrs}

Let $\theta=(\beta^\top,\gamma^\top)^\top$ be the unknown parameter vector
that indexes the extreme value regression model (\ref{extreme value-max})-(\ref{linkphi}). 
In what follows, $\nu$ represents the scalar parameter of interest  
and $\psi=(\psi_1,\ldots,\psi_s)^\top$ is the nuisance parameter; note that $1+s=k+m$. 
Let $\ell(\theta)$ be the log-likelihood function and let $I$ and $J$ be the
expected and observed information matrices, respectively, which are presented in \citet[p. 584]{FERRARI2012}.
Furthermore, let $J_{\psi \psi}$ denote the $s \times s$
observed information matrix corresponding to $\psi$. Similarly, 
$A_{\psi \psi}$ denotes a matrix formed from the $(r+s) \times (r+s)$ matrix 
$A$ by dropping the rows and columns that correspond to the parameter of interest. 
Additionally, hat and tilde symbols indicate evaluation at the unrestricted 
($\widehat\theta$) and at the restricted ($\widetilde\theta$) maximum likelihood 
estimator of $\theta$ for a given $\psi$, respectively. For instance,
$\widehat{I}=I(\widehat{\theta})$, $\widetilde{I}=I(\widetilde{\theta})$, and
$\widehat{J}=J(\widehat{\theta})$. 

Inference on $\nu$ can be based on the signed likelihood ratio statistic
\[
R=R(\nu) = {\rm sgn}(\widehat{\nu} - \nu)\bigl\{2\bigl[\ell(\widehat{\theta}) -
\ell(\widetilde{\theta})\bigr]\bigr\}^{1/2}.
\]
It is well known that $R$ has the first order of
accuracy only, i.e.~$R$ is approximately distributed as a standard normal distribution with
error of order $O(n^{-1/2})$. This approximation can be inaccurate when the sample
size is small. The accuracy of this approximation can be improved by making suitable adjustments to $R$.

\citet{BARNDORFF-NIELSEN1986,BARNDORFF-NIELSEN1991} proposed a modified signed likelihood ratio statistic, $R^*$,
which has an asymptotic normal distribution with approximation error of order $O(n^{-3/2})$.
Significance tests and confidence limits based on the normal approximation to the distribution of $R^*$ are, in general, extremely accurate even for very small sample sizes. 
The proposed modification is defined in terms of derivatives of the log-likelihood function with respect to the sample space, 
more precisely, derivatives with respect to the maximum likelihood estimator of $\theta$. To obtain the derivatives, one needs to specify an ancillary statistic, for example ${\mathbf a}$, such that $(\widehat \theta, {\mathbf a})$ is a sufficient statistic. 
Thus, the log-likelihood function can be written as $\ell(\theta;\widehat \theta,{\mathbf a})$ and one can derive 
$\ell(\theta;\widehat \theta,{\mathbf a})$ with respect to $\widehat \theta$. 
In the following, we use the standard notation for derivatives of the log-likelihood function:
$
\ell_{\theta}(\theta)
=\ell_{\theta}(\theta;\widehat{\theta},{\mathbf a})
=\partial \ell(\theta;\widehat{\theta},{\mathbf a}) / \partial \theta^{\!\top}, 
$
and
$
\ell_{;\widehat{\theta}}(\theta)
=\ell_{;\widehat{\theta}}(\theta;\widehat{\theta},{\mathbf a})
=\partial \ell(\theta;\widehat{\theta},{\mathbf a}) / \partial \widehat{\theta}. 
$
Hence, $\ell_{\theta}(\theta)$ is a $(k+m)$  column vector, $\ell_{;\widehat{\theta}}(\theta)$ is a $(k+m)$ row vector,
and  $\ell_{\theta;\theta}(\theta)$ and $\ell_{\theta;\widehat{\theta}}(\theta)$ are $(k+m) \times (k+m)$ matrices.

Barndorff-Nielsen's modified signed likelihood ratio statistics is given by
\begin{equation}\label{barndorff}
R^*=R+\frac{1}{R}\ln \left( \left|\frac{U}{R}\right| \right),
\end{equation}
where
$$ 
U=
\frac{
\left|\begin{array}{c c}
\ell_{;\widehat{\theta}}(\widehat{\theta})-\ell_{;\widehat{\theta}}(\widetilde{\theta})
\\
\ell_{\psi;\widehat{\theta}}(\widetilde{\theta})
\end{array}
\right|}
{|\widetilde J_{\psi \psi}|^{1/2}|\widehat{J}|^{1/2}}.
$$
Here $\ell_{\psi;\widehat{\theta}}(\widetilde{\theta})$ is obtained from $\ell_{\theta;\widehat{\theta}}(\widetilde{\theta})$ excluding 
the row that corresponds to the parameter of interest $\nu$.
For the one-parameter case, i.e. the case in which there is no nuisance parameter, we have 
$$
U=|\ell_{;\widehat{\theta}}(\widehat{\theta})-\ell_{;\widehat{\theta}}(\theta)|\;|\widehat{J}|^{-1/2}.
$$

Because the computation of $U$ involves sample space derivatives, which, in turn, requires the determination of a suitable ancillary
statistic, it can be difficult or even impossible to obtain. Most of the applications of Barndorff-Nielsen's adjustment concern exponential family models and transformation models (e.g. \citealp{WONG2002}, \citealp{LARSEN2002}, \citealp{MELO2010}
and \citealp{CORTESE2013}). Approximations to $R^*$ that avoid the need for sample space derivatives have been proposed by several authors. The different approximations lead to alternative adjusted signed likelihood ratio statistics. It is our aim to obtain and compare the different adjusted statistics for the extreme value regression models introduced in Section~\ref{model}.

\citet{DICICCIO} proposed an approximation to Barndorff-Nielsen's modified statistic that requires that the parameter of interest and the nuisance parameters be globally orthogonal, i.e., the Fisher information matrix is block-diagonal.
It should be noted that it is always possible to obtain an orthogonal parameterization when the parameter of interest is scalar
\citep{COX1987}. The adjusted signed likelihood ratio statistic proposed by \citet{DICICCIO}, 
denoted here by $R^*_0$, is defined by (\ref{barndorff}) with $U$ replaced by
\begin{equation}
\label{U0}
U_0=
{\ell}_\nu(\widetilde \theta)
\Biggl (
\frac{|\widetilde J_{\psi \psi}||\widehat I_{\nu\nu}|}
{|\widehat J||\widetilde I_{\nu\nu}|}
\Biggr )^{1/2},
\end{equation}
the sign of $U_0$ being taken to be the same as that of $R$, where $\ell_\nu(\theta)$ is the derivative of the log-likelihood function with respect to the parameter of interest. According to \citet[Section 7.5.2]{SEVERINI2000}, 
$R_0^*=R^*+O_p(n^{-1})$ for $\nu$ of the form $\widehat\nu=\nu+O(n^{-1})$.

The adjusted signed likelihood ratio statistic proposed by \citet{SKOVGAARD1996}, denoted here by $\overline{R}^*$, is based on 
covariances and is given by
(\ref{barndorff}) with $U$ replaced by
\begin{equation}\label{overlineU}
\overline{U}
=
\left|
\begin{pmatrix}
q^{\!\top} \\
\Upsilon_{\psi}
\end{pmatrix}
\right|
|\widehat I|^{-1}
|\widehat J|^{1/2}
|\widetilde J_{\psi \psi}|^{-1/2},
\end{equation}
where $q=q(\widehat \omega, \widetilde \omega)$ and $\Upsilon=\Upsilon(\widehat \omega, \widetilde \omega)$, with
$$ 
q(\omega_1,\omega)={\rm E}_{\omega_1}[U(\omega_1) \ (\ell(\omega_1)-\ell(\omega))], 
\qquad
\Upsilon(\omega_1,\omega)={\rm E}_{\omega_1}[U(\omega_1) \ U^{\top}(\omega)].
$$
Also, $\Upsilon_{\psi}$ is the matrix $\Upsilon$ without the row that corresponds to the parameter of interest $\nu$.
According to \citet[Section 7.5.4]{SEVERINI2000}, 
$\overline{R}^*=R^*+O_p(n^{-1})$ for $\nu$ of the form $\widehat\nu=\nu+O(n^{-1})$.

The adjusted signed likelihood ratio statistic proposed by \citet{SEVERINI1999}, denoted here by $\widehat{R}^*$, is based on 
empirical covariances and is given by (\ref{barndorff}) with $U$ replaced by
\begin{equation}\label{widehatU}
\widehat{U}
=
\left|
\begin{pmatrix}
\widehat{q}^{\!\top} \\
\widehat{\Upsilon}_{\psi}
\end{pmatrix}
\right|
|\widehat I|^{-1}
|\widehat J|^{1/2}
|\widetilde J_{\psi \psi}|^{-1/2},
\end{equation}
where $\widehat q$ and $\widehat \Upsilon$ are empirical estimates of  $q$ and $\Upsilon$ given by
\begin{eqnarray}
\widehat q=\sum_{t=1}^n(\ell_t(\widehat\theta)-\ell_t(\widetilde\theta))\frac{\partial \ell_t(\widehat\theta)}{\partial \theta ^{\!\top}}
\quad
{\rm and}
\quad
\widehat \Upsilon = \sum_{t=1}^n  \frac{\partial \ell_t(\widetilde \theta)}{\partial \theta} \frac{\partial \ell_t(\widehat\theta)}{\partial \theta ^{\!\top}}.
\end{eqnarray}
According to \citet[Section 7.5.5]{SEVERINI2000}, 
$\widehat{R}^*=R^*+O_p(n^{-1})$ for $\nu$ of the form $\widehat\nu=\nu+O(n^{-1})$. 

The adjusted signed likelihood ratio statistic given by Fraser, Reid \& Wu \citep{FRASER1999}, denoted here by $\widetilde{R}^*$, 
is based on an approximate ancillary statistic, and is defined by (\ref{barndorff}) with $U$ replaced by
$$
\widetilde{U}=
\frac{
\left|\begin{array}{c c}
\widetilde \ell_{;\widehat{\theta}}(\widehat{\theta})- \widetilde \ell_{;\widehat{\theta}}(\widetilde{\theta})
\\
\widetilde \ell_{\psi;\widehat{\theta}}(\widetilde{\theta})
\end{array}
\right|}
{|\widetilde J_{\psi \psi}|^{1/2}|\widehat{J}|^{1/2}},
$$
where
$$
\widetilde \ell_{;\widehat{\theta}}(\widehat{\theta})- \widetilde \ell_{;\widehat{\theta}}(\widetilde{\theta})=
\left( \ell_{;y}(\widehat\theta) -  \ell_{;y}(\widetilde \theta) \right)\widehat{V}\left( \ell_{\theta;y}(\widehat \theta)\widehat{V} \right)^{-1}\widehat J 
$$
is an $1 \times (k+m)$ vector,
$$
\widetilde \ell_{\theta;\widehat \theta}(\widetilde \theta)
=\ell_{\theta;y}(\widetilde \theta)\widehat{V}\left( \ell_{\theta;y}(\widehat \theta)\widehat{V} \right)^{-1}\widehat J
$$
is a $(k+m) \times (k+m)$ matrix, $\ell_{;y}$ is the derivative of the log-likelihood function with respect to $y$, $\ell_{\theta;y}$ is its second derivative with respect to $\theta$ and $y$, and $\widehat V$ is an $n \times (k+m)$ matrix with the element $(i,j)$ given by 
$\left( -\left( \partial F(y_i,\theta) / \partial \theta_j \right) / f(y_i,\theta) \right)\bigl|_{\theta = \widehat\theta}$, $F$ is the cumulative distribution function of the model and $f$ is the corresponding probability density function.
According to \citet[Section 7.5.3]{SEVERINI2000}, 
$\widetilde{R}^*=R^*+O_p(n^{-3/2})$ for $\nu$ of the form $\widehat\nu=\nu+O(n^{-1})$. 
It should be noted that the approximation error for $\widetilde{R}^*$ is smaller than those corresponding to $R^*_0$, $\overline R^*$ and $\widehat R^*$, which are all $O_p(n^{-1})$.

\section{Adjusted signed likelihood ratio statistics for extreme value regression models}\label{aslrsevrm}

First, we note that the adjusted signed likelihood ratio statistics presented in the previous section require the computation
of the observed and/or expected information matrices. Both matrices for extreme value regression models can be found in \citet{FERRARI2012}. 

It is always possible to obtain Barndorff-Nielsen's modified statistic $R^*$ in homoskedastic linear regression models by considering the ancillary statistic ${\mathbf a}=(a_1,\ldots,a_n)$, where $a_t=(y_t-x_t^\top\widehat\beta)/\widehat{\sigma}$.
The linear maximum extreme value regression model with constant dispersion is a special case of model (\ref{extreme value-max})-(\ref{linkphi}), in where $g$ and $h$ are the identity function, $\eta_t=x_t^{\!\top}\beta$ and $\delta_t=\sigma$. For this model, 
$$
\ell_{;\widehat{\theta}}(\widehat{\theta})-\ell_{;\widehat{\theta}}(\widetilde{\theta})
=
\iota^{\!\top} \left( \widehat{\sigma}^{-1} (-{\mathcal I} +\widehat{\breve {\mathcal Z}}) -  \widetilde{\sigma}^{-1} (-{\mathcal I} +\widetilde{\breve {\mathcal Z}})\right)
\begin{pmatrix}
\widehat X & \widehat{\mathcal Z}\iota 
\end{pmatrix},
$$
$$
\ell_{ \theta;\widehat{\theta}}(\widetilde \theta)
=\frac{\partial^2 \ell(\theta;\widehat\theta,a) }{\partial \theta \partial \widehat \theta}\biggr|_{\widetilde \theta}
=
\widetilde A
\begin{pmatrix}
\widehat X & \widehat{\mathcal Z} \iota 
\end{pmatrix}, 
$$
where
$$
A=
\sigma^{-2}
\begin{pmatrix}
 X^{\!\top}  {\breve {\mathcal Z}} 
\\
 \iota^{\!\top} ({\mathcal I}-{\breve {\mathcal Z}} + {\mathcal Z} {\breve {\mathcal Z}})
\end{pmatrix},
$$
$X=(x_1^\top,\ldots,x_n^\top)^\top$,
$\iota$ is an $n$-dimensional column vector of ones,
${\mathcal I}$ is the $n \times n$ identity matrix, 
${\mathcal Z}={\rm diag} \{z _1,\ldots,z_n\}$,
${\breve {\mathcal Z}} = {\rm diag}\{\exp(-z_1),\ldots,\exp (-z_n)\}$
and $z_t = (y_t-\mu _t)/\sigma _t$.
Hence, Barndorff-Nielsen's modified statistic $R^*$ can be written as in (\ref{barndorff}) with
$$
U 
=
\left| 
\begin{pmatrix}
\iota^{\!\top} \left(\widehat \sigma^{-1}(-{\mathcal I}+\widehat {\breve {\mathcal Z}}) - \widetilde \sigma^{-1}(-{\mathcal I}+\widetilde {\breve {\mathcal Z}}) \right)\\ 
\widetilde A_{\psi} 
\end{pmatrix}
\begin{pmatrix}
\widehat X & \widehat{\mathcal Z} \iota
\end{pmatrix} 
\right|  |\widetilde J_{\psi \psi}|^{-1/2} |\widehat J|^{-1/2},
$$   
where $A_\psi$ is obtained from $A$ by excluding the row that corresponds to $\nu$. 

We now obtain the aforementioned approximations to $R^*$ for the maximum extreme value regression model (\ref{extreme value-max})-(\ref{linkphi}). It should be noted that, hereafter the linearity and constant dispersion assumptions are no longer required.

One should take the following steps to obtain DiCiccio \& Martin's adjusted signed likelihood ratio statistic $R_0^*$ \citep{DICICCIO}: find a parameterization for which the parameter of interest is orthogonal to the 
remaining parameters; re-write (\ref{extreme value-max})-(\ref{linkphi}) in the orthogonal parameterization;
obtain the log-likelihood function, its derivative with respect to the parameter of interest,
the observed information matrix and the element of the expected information matrix that corresponds to the parameter of interest;
compute (\ref{U0}). It should be noted that the adjustment term depends on the choice of the parameter of interest and 
the systematic sub-models. In the Appendix we present the steps above for inference on one of the components of the $\beta$ vector
in a model with linear specification for the location and dispersion sub-models.

Skovgaard's adjusted signed likelihood ratio statistic $\overline{R}^*$ \citep{SKOVGAARD1996} is obtained by replacing $q$ and 
$\Upsilon$ in (\ref{overlineU})  
with
$$
\overline q=
\left[
\begin{array}{c}
\widehat X^{\!\top} \widehat \Phi^{-1} \widehat T C \left( {\mathcal I}-M\breve D \right)\iota \\
\widehat Z^{\!\top} \widehat \Phi^{-1} \widehat H   \left( C({\mathcal E}{\mathcal I}+N\breve D)-{\mathcal I} \right)\iota
\end{array}
\right]
$$
and
$$
\label{eq:UpsilonEV}
\overline \Upsilon=
\left[
\begin{array}{c c}
\widehat X^{\!\top} \widehat \Phi^{-1} \widehat T CM\breve D \widetilde T \widetilde \Phi^{-1}\widetilde X   &
\widehat X^{\!\top} \widehat \Phi^{-1} \widehat T C\{ {\mathcal I} + \breve D  (-M - CN + MD) \}\widetilde H \widetilde \Phi^{-1} \widetilde Z \\
-\widehat Z^{\!\top} \widehat \Phi^{-1} \widehat H CN\breve D \widetilde T \widetilde \Phi^{-1} \widetilde X     \quad&
\widehat Z^{\!\top} \widehat \Phi^{-1} \widehat H C \{ {\mathcal E}{\mathcal I} + \breve D (N + CP - ND) \} \widetilde H \widetilde \Phi^{-1} \widetilde Z 
\end{array}
\right],
$$
respectively, where 
$T = {\rm diag}\{ 1/g'(\mu_{1}), \ldots, 1/g'(\mu_{n})\}$, 
$H = {\rm diag}\{ 1/h'(\sigma_{1}), \ldots, 1/h'(\sigma_{n})\}$, 
$\Phi={\rm diag}\{ \sigma_{1}, \ldots, \sigma_{n}\}$,
$C={\rm diag}\{\widehat \sigma_{1}/\widetilde \sigma_1,\ldots,\widehat \sigma_{n}/\widetilde \sigma_n\}$, 
$D={\rm diag}\{(\widehat \mu_{1}-\widetilde \mu_1)/\widetilde \sigma_1,\ldots,(\widehat \mu_{n}-\widetilde \mu_n)/\widetilde \sigma_n\}$,
$\breve D={\rm diag}\{\exp(-(\widehat \mu_{1}-\widetilde \mu_1)/\widetilde \sigma_1),\ldots,\exp(-(\widehat \mu_{n}-\widetilde \mu_n)/\widetilde \sigma_n)\}$, 
$M={\rm diag}\{\Gamma(1+\widehat \sigma_{1}/\widetilde \sigma_1),\ldots,\Gamma(1+\widehat \sigma_{n}/\widetilde \sigma_n)\}$,
$N={\rm diag}\{\Gamma^{(1)}(1+\widehat \sigma_{1}/\widetilde \sigma_1),\ldots,\Gamma^{(1)}(1+\widehat \sigma_{n}/\widetilde \sigma_n)\}$,
$P={\rm diag}\{\Gamma^{(2)}(1+\widehat \sigma_{1}/\widetilde \sigma_1),\ldots,\Gamma^{(2)}(1+\widehat \sigma_{n}/\widetilde \sigma_n)\}$,
and the other matrices are as defined above. Here, $\Gamma(\cdot)$ is the gamma function, and $\Gamma^{(1)}(\cdot)$ and $\Gamma^{(2)}(\cdot)$
denote the first and second derivatives of the function, respectively.

To obtain Severini's adjusted signed likelihood ratio statistic $\widehat{R}^*$ \citep{SEVERINI1999}, $\widehat q$ and $\widehat \Upsilon$ are replaced in (\ref{widehatU}) by
\begin{eqnarray*}
\widehat q =
\begin{pmatrix}
 \widehat X \widehat \Phi^{-1} \widehat T\left( \widehat {\mathfrak L} - \widetilde {\mathfrak L} \right)({\mathcal I} - \widehat {\breve {\mathcal Z}})\iota\\
 \widehat Z \widehat \Phi^{-1} \widehat H\left( \widehat {\mathfrak L} - \widetilde {\mathfrak L} \right)(-{\mathcal I} + \widehat {\mathcal Z}-\widehat {\mathcal Z}\widehat {\breve {\mathcal Z}})\iota 
\end{pmatrix}
\quad
{\rm and}
\quad\widehat \Upsilon =\left[
\begin{array}{c c}
\widehat \Upsilon_{\beta \beta } & \widehat \Upsilon_{\beta \gamma }\\
\widehat \Upsilon_{\gamma \beta } & \widehat \Upsilon_{\gamma \gamma }
\end{array}
\right],
\end{eqnarray*}
where
$$
\widehat \Upsilon_{\beta \beta } =
\widetilde X^{\!\top}\widetilde \Phi ^{-1}\widetilde T({\mathcal I}-\widetilde {\breve {\mathcal Z}}) \ 
({\mathcal I}-\widehat {\breve {\mathcal Z}})\widehat \Phi ^{-1} \widehat T \widehat X,
$$
$$
\widehat \Upsilon_{\beta \gamma }=
\widetilde X^{\!\top}\widetilde \Phi ^{-1}\widetilde T({\mathcal I}-\widetilde {\breve {\mathcal Z}}) \
(-{\mathcal I}+\widehat {\mathcal Z}-\widehat {\mathcal Z}\widehat {\breve {\mathcal Z}})\widehat \Phi ^{-1} \widehat H \widehat Z ,
$$
$$
\widehat \Upsilon_{\gamma \beta }=
\widetilde Z \widetilde H \widetilde \Phi ^{-1}(-{\mathcal I}+\widetilde {\mathcal Z}-\widetilde {\mathcal Z}\widetilde {\breve {\mathcal Z}}) \
({\mathcal I}-\widehat {\breve {\mathcal Z}})\widehat \Phi ^{-1} \widehat T \widehat X,
$$
$$
\widehat \Upsilon_{\gamma \gamma }=
\widetilde Z \widetilde H \widetilde \Phi ^{-1}(-{\mathcal I}+\widetilde {\mathcal Z}-\widetilde {\mathcal Z}\widetilde {\breve {\mathcal Z}})
(-{\mathcal I}+\widehat {\mathcal Z}-\widehat {\mathcal Z}\widehat {\breve {\mathcal Z}})\widehat \Phi ^{-1} \widehat H \widehat Z, 
$$
${\mathfrak L}= {\rm diag} \{ \ell_1,\ldots,\ell_n \}$, with
$$\ell_t=\ell_t(\mu_t,\sigma_t)=\ln f(y_t;\mu_t,\sigma_t)=-\ln \sigma_t -(y_t-\mu_t)/(\sigma_t)-\exp \left( -(y_t-\mu_t)/(\sigma_t) \right),$$  
and the other matrices are as defined above.

Finally, Fraser, Reid \& Wu's adjusted signed likelihood ratio statistic $\widetilde{R}^*$ \citep{FRASER1999} is given by (\ref{barndorff}) with $U$ replaced by
$$
\widetilde U 
=
\left| 
\begin{pmatrix}
\iota \left(\widehat \Phi^{-1}(-{\mathcal I}+\widehat {\breve{\mathcal Z}}) - \widetilde \Phi^{-1}(-{\mathcal I}+\widetilde {\breve{\mathcal Z}}) \right)\\ 
\widetilde A_{\psi} 
\end{pmatrix}
\begin{pmatrix} \widehat T \widehat X & \widehat{\mathcal Z} \widehat H  \widehat Z \end{pmatrix}
\begin{pmatrix}
\widehat A
\begin{pmatrix} \widehat T \widehat X & \widehat{\mathcal Z} \widehat H  \widehat Z \end{pmatrix}
\end{pmatrix}^{-1}
   \right|  |\widetilde J_{\psi \psi}|^{-1/2} |\widehat J|^{1/2},
$$  
where
$$
A=\begin{pmatrix}
X^{\!\top}  T  \Phi^{-2}  {\breve{\mathcal Z}}\\
 Z^{\!\top}  H  \Phi^{-2}({\mathcal I}- {\breve{\mathcal Z}}+ {\mathcal Z} {\breve{\mathcal Z}})
\end{pmatrix},
$$
$\widehat A=A(\widehat\theta)$, $\widetilde A=A(\widetilde \theta)$, $\widetilde \ell_{\psi;\widehat{\theta}}(\widetilde{\theta})$ is obtained from  $\widetilde \ell_{\theta;\widehat{\theta}}(\widetilde{\theta})$ by excluding the row corresponding the parameter of 
interest $\nu$, and the other matrices are as  defined above.

All of the formulas presented above are valid for maximum extreme value regression models. If, instead, the data are assumed to follow a  
minimum extreme value regression model, the formulas can be easily adapted. Recall that, from (\ref{min-max}), 
if $y_t \sim EV_{min}(\mu_t,\sigma_t)$, then $y_t^*=-y_t$ has a maximum extreme value distribution with parameters $-\mu_t$ and 
$\sigma_t$ ($-y_t \sim VE_{max}(-\mu_t,\sigma_t)$). Hence, one can obtain all of the adjusted statistics by taking the response variable 
$y_t^*$, and the systematic components $g^*(\mu_t)=\mu_t^*=\eta^*_t=-g^{-1}(\eta(x_t,\beta))$ and $h^*(\sigma_t)=h(\sigma_t)=\delta(z_t,\gamma)$.
In other words, adjusted signed likelihood ratio statistics for the minimum extreme value regression model (\ref{extreme value-min}) with systematic 
components (\ref{linkmu}) and (\ref{linkphi}) can be computed from the results above by changing the signs of the observations on the response variable 
and using an identity link function for the location parameter with the modified predictor $\eta^*_t$. 
We note that Fraser, Reid \& Wu's adjusted statistic
for the minimum extreme value distribution (iid case) under censoring was obtained by \citet{WONG2000}.

\section{Monte Carlo simulation results}\label{MonteCarlo}

We now present a Monte Carlo simulation study to evaluate the finite-sample performance of  
the signed likelihood ratio test ($R$), Barndorff-Nielsen's modified version ($R^*$) and the tests that use approximations to $R^*$, namely
DiCiccio \& Martin's $R^*_0$, 
Skovgaard's $\overline{R}^*$, 
Severini's $\widehat{R}^*$ and 
Fraser, Reid \& Wu's $\widetilde{R}^*$, 
in extreme value regression models.
The maximum likelihood estimates of the parameters are obtained by numerically maximizing the log-likelihood function using the nonlinear 
quasi-Newton BFGS method with analytical derivatives (see, for instance, 
\citealp{PRESS}), implemented in the function MaxBFGS in the matrix programming language {\tt Ox} 
\citep{DOORNIK}. The number of Monte Carlo replicates is 10,000 and the nominal levels of the tests are $\alpha = 10\%$, 5\% and 1\%.

At the outset, we consider the maximum extreme value regression model with constant dispersion
and location sub-model
$$
\mu_t = \beta_0 + \beta_1 x_{t1} + \beta_2 x_{t2}
                                    + \beta_3 x_{t3} + \beta_4 x_{t4};
$$
`model 1'. The hypothesis being tested is ${\mathcal H}_0: \beta_1 \leq 0$ and the alternative hypothesis is ${\mathcal H}_1: \beta_1 > 0$.
We set $\beta_0=1$, $\beta_1=0$, $\beta_2=1$, $\beta_3=6$, $\beta_4=-3$ and $\sigma=1$.
The values of the covariates are randomly drawn from the 
${\mathcal U}(-0.5, \ 0.5)$ distribution. The sample sizes range from 15 to 200.

Table~\ref{tab:sinalizada-linear-homo_tamanho} presents the null rejection rates of the six tests. 
The figures in this table indicate that the signed likelihood ratio test is liberal, more so when the sample size is small.
It is noteworthy that all of the adjusted tests are much less size distorted.
For instance, when $n=20$ and $\alpha=5\%$, the null rejection rate is $8.1\%$ whereas the rejection rates of the adjusted tests are
$5.0\%$  (Barndorff-Nielsen's $R^*$, DiCiccio \& Martin's $R^*_0$  and	Fraser, Reid \& Wu's $\widetilde{R}^*$),  $4.0\%$ (Skovgaard's $\overline{R}^*$) 
and $3.8\%$ (Severini's $\widehat{R}^*$). As expected, the null
rejection rates of all the tests are near the nominal levels for large samples.

\begin{figure}[!ht]
  \centering
  \includegraphics[height=70mm,width=165mm]{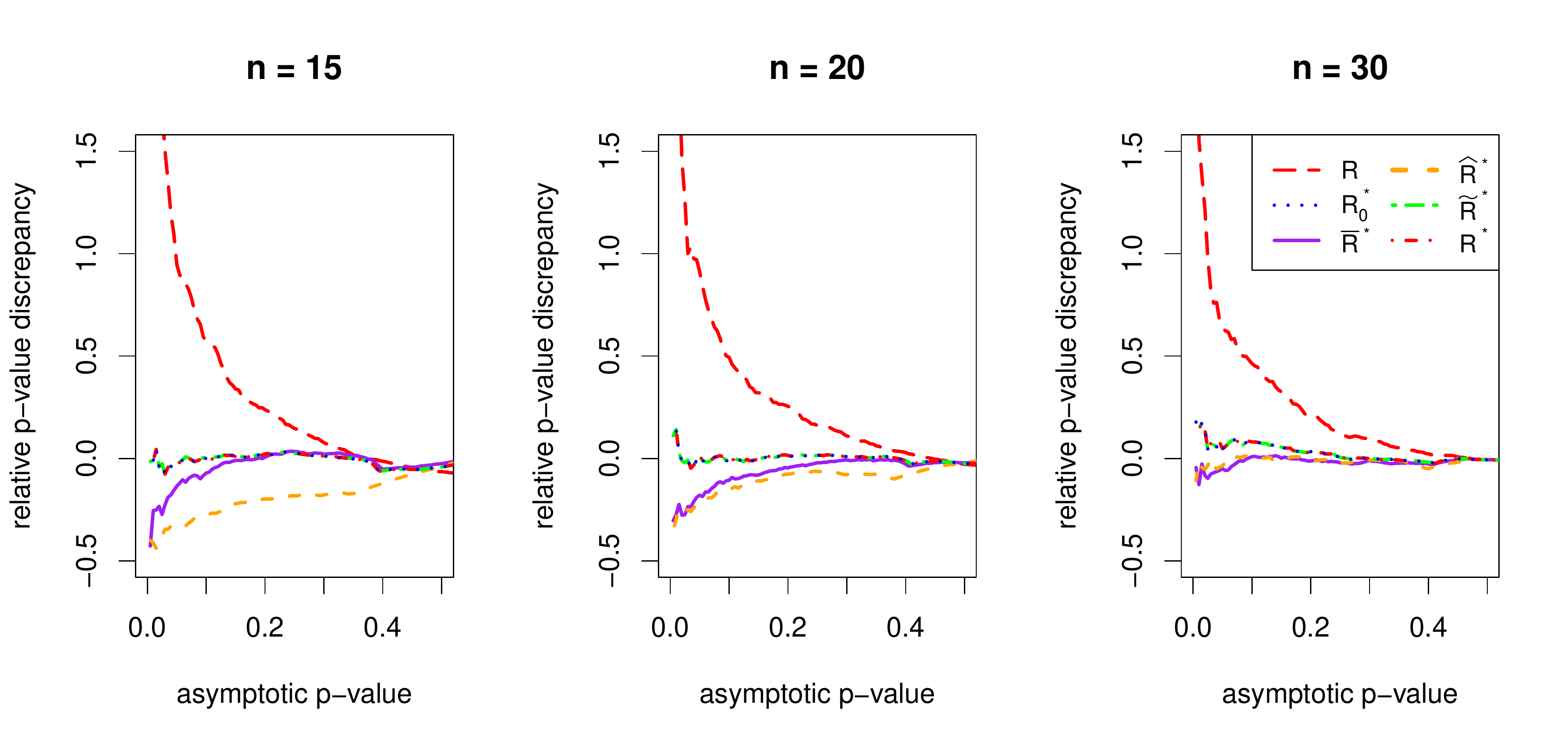}
  \caption{Relative p-value discrepancy plots, model 1.}
  \label{fig:sinalizada-linear-homo-valorp}
\end{figure}

Figure~\ref{fig:sinalizada-linear-homo-valorp} presents plots of relative p-value discrepancies versus asymptotic p-values for
small to moderate sample sizes, namely $n=15, \ 20, \ 30$. 
The relative p-value discrepancy is the difference between the exact and the asymptotic p-values divided by the latter.
It should be noted that, on the left-hand side of the plots (asymptotic p-value smaller than 
$10\%$) the relative discrepancy line of the signed likelihood ratio test is well above the zero line.
Hence, the plots show that the signed likelihood ratio test is anti-conservative, as anticipated by the figures presented in 
Table~\ref{tab:sinalizada-linear-homo_tamanho}. 
The adjusted tests are clearly better behaved, particularly Barndorff-Nielsen's $R^*$, DiCiccio \& Martin's $R^*_0$ and Fraser, Reid and Wu's $\widetilde{R}^*$. 
Their relative p-value discrepancies are not far from zero, even for very small asymptotic p-values (e.g., 0.1\%) and small sample size (e.g., $n=15$).
The plots also show that the other tests, Skovgaard's $\overline{R}^*$ and Severini's $\widehat{R}^*$, are conservative when $n=15$ and $20$. The conservative 
behavior is more pronounced for $\overline{R}^*$.  When $n=30$ the p-value discrepancies for both tests are very small.

We now turn to power simulations. As noted above, the different tests have different sizes when the asymptotic normal distribution is used as an approximation to the true null distribution of the statistics. To guarantee that all tests have the correct size, we first simulate 500,000 samples under the null hypothesis being tested. We then obtain the empirical distribution of each statistic under the null hypothesis, and we estimate the exact critical values of tests that give the correct size for the chosen nominal level. The exact critical values are used to evaluate the power of the tests. This strategy is used for all of the power simulations in this paper. Table~\ref{tab:sinalizada-linear-homo_poder} presents the rejection rates of the tests under the alternative hypotheses $\beta_2=\epsilon$ for different values of $\epsilon>0$.  
It is apparent that all of the tests have very similar powers.
\begin{table}[htp]
\centering
\renewcommand{\arraystretch}{1.1}
\caption{Null rejection rates (\%), model 1}\label{tab:sinalizada-linear-homo_tamanho}
\scriptsize
\begin{tabular}
{c| c c c c c c| c c c c c c|  c c c c c c}
\hline \hline
& \multicolumn{6}{c|}{$\alpha = 10\%$}
& \multicolumn{6}{c|}{$\alpha = 5\%$}
& \multicolumn{6}{c}{$\alpha = 1\%$}
\\
\cline{2-19}
& \multicolumn{6}{c|}{}
& \multicolumn{6}{c|}{}
& \multicolumn{6}{c}{}
\\[-0.1in]
$n$ & $R$	&	$R^*_0$ &	$\overline{R}^*$	&	$\widehat{R}^*$	&	$\widetilde{R}^*$ &	$R^*$	
    &	$R$	&	$R^*_0$ &	$\overline{R}^*$	&	$\widehat{R}^*$	&	$\widetilde{R}^*$ &	$R^*$	
    &	$R$	&	$R^*_0$ &	$\overline{R}^*$	& $\widehat{R}^*$	&	$\widetilde{R}^*$ &	$R^*$	
            \\
\hline            
15	&	13.8	&	10.0	&	9.0	  &	6.8	  &	10.0	&	10.0	&	8.4	&	4.9	&	4.2	&	3.3	&	4.9	&	4.9	&	2.7	&	1.0	&	0.7	&	0.6	&	1.0	&	1.0	\\
20	&	13.5	&	9.9	  &	8.8	  &	8.3	  &	9.9	  &	9.9	  &	8.1	&	5.0	&	4.0	&	3.8	&	5.0	&	5.0	&	2.3	&	1.1	&	0.8	&	0.7	&	1.1	&	1.1	\\
30	&	13.7	&	10.8	&	10.1	&	10.1	&	10.8	&	10.8	&	7.7	&	5.2	&	4.8	&	4.8	&	5.2	&	5.2	&	2.2	&	1.2	&	0.9	&	1.0	&	1.2	&	1.2	\\
40	&	12.2	&	9.6	  &	8.8	  &	9.7	  &	9.6	  &	9.6	  &	6.7	&	4.8	&	4.3	&	4.8	&	4.8	&	4.8	&	1.7	&	1.2	&	1.0	&	1.1	&	1.2	&	1.2	\\
100	&	10.5	&	9.6  	&	9.2	  &	9.5	  &	9.6	  &	9.6	  &	5.2	&	4.7	&	4.6	&	4.8	&	4.7	&	4.7	&	1.2	&	0.9	&	0.8	&	0.9	&	0.9	&	0.9	\\
200	&	10.5	&	10.1	&	10.0	&	10.2	&	10.1	&	10.1	&	5.3	&	5.0	&	4.8	&	5.0	&	5.0	&	5.0	&	1.2	&	1.0	&	0.9	&	1.1	&	1.0	&	1.0	\\
\hline \hline
\end{tabular}
\end{table}
\begin{table}[!ht]
\begin{center}
\caption{Non-null rejection rates (\%), model 1, $n=20$, $\alpha=10\%$} \label{tab:sinalizada-linear-homo_poder}
\small
\begin{tabular}
{c| c c c c c c c c c c c c}
\hline \hline
$\epsilon$	&	0.1	&	0.5	&	1.0	&	2.0	&	3.0	\\
\hline
$R$	&	11.7	&	22.3	&	39.8	&	76.5	&	94.6	\\
$R^*_0$	&	11.7	&	22.2	&	39.8	&	76.9	&	94.8	\\
$\overline{R}^*$	&	11.7	&	22.1	&	39.8	&	76.8	&	94.8	\\
$\widehat{R}^*$	&	11.7	&	21.6	&	37.9	&	72.5	&	91.8	\\
$\widetilde{R}^*$ 	&	11.7	&	22.2	&	39.8	&	76.9	&	94.8	\\
$R^*$	&	11.7	&	22.2	&	39.8	&	76.9	&	94.8	\\
\hline \hline
\end{tabular}
\end{center}
\end{table}

\begin{figure}[!ht]
  \centering
  \includegraphics[height=70mm,width=165mm]{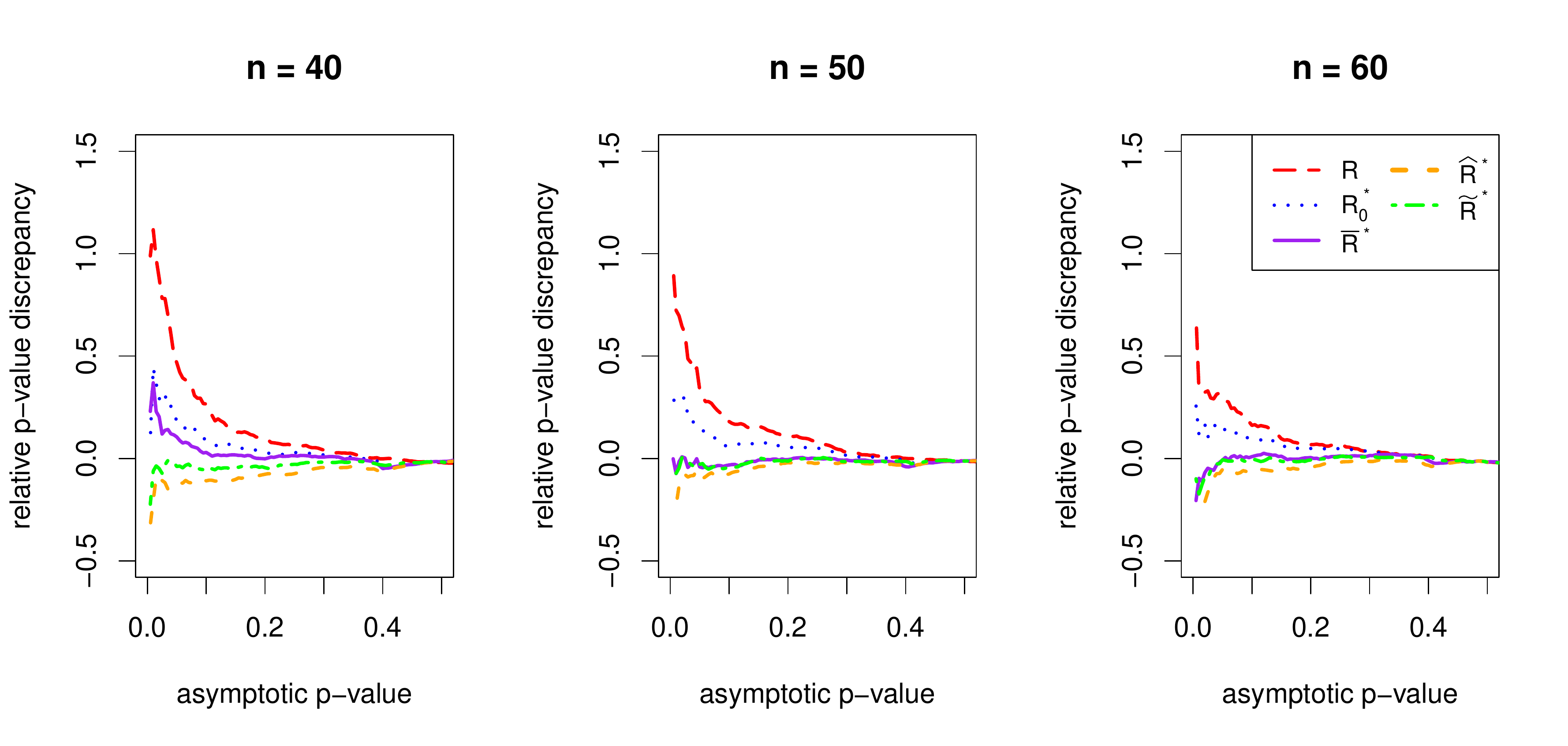}
  \caption{Relative p-value discrepancy plots, model 2.}
  \label{fig:sinalizada-linear-hetero-valorp}
\end{figure}

We now consider the maximum extreme value regression model with location and dispersion sub-models
$$
\mu_t = \beta_0 + \beta_1 x_{t1} + \beta_2 x_{t2} + \beta_3 x_{t3}
$$
and
$$
\ln (\sigma_t) = \gamma_0 + \gamma_1 z_{t1} + \gamma_2 z_{t2};
$$
`model 2'. 
The null and alternative hypotheses are ${\mathcal H}_0: \beta_3 \leq 0$ and ${\mathcal H}_1: \beta_3 > 0$ .
Here, we set $\beta_0=1$, $\beta_1=1$, $\beta_2=6$, $\beta_3=0$, $\gamma_0=1$ and $\gamma_1=\gamma_2=0.1$,
where the values of the covariates are drawn from a ${\mathcal U}(-0.5, \ 0.5)$ distribution, and the sample sizes are 40, 50, 60 and 70.
For this model and for model 3 (below), it is not possible to obtain the adjusted signed likelihood ratio statistic $R^*$. 
Therefore, we only consider the approximations to $R^*$ presented in the previous section.
The size and power simulation results are summarized in Tables~\ref{tab:sinalizada-linear-hetero-H0beta_tamanho} and
\ref{tab:sinalizada-linear-hetero-H0beta_poder} and in Figure~\ref{fig:sinalizada-linear-hetero-valorp}. Again, it is clear that
the signed likelihood ratio test is oversized when the sample is not large, the proposed adjustments are effective in pushing 
the size to the chosen nominal level, all of the rejection rates converge to the nominal level as $n$ grows, and the tests are equally powerful.  
\begin{table}[htp]
\centering
\renewcommand{\arraystretch}{1.1}
\caption{Null rejection rates (\%), model 2}\label{tab:sinalizada-linear-hetero-H0beta_tamanho}
\small
\begin{tabular}
{c| c c c c c | c c c c c|  c c c c c}
\hline \hline
& \multicolumn{5}{c|}{$\alpha = 10\%$}
& \multicolumn{5}{c|}{$\alpha = 5\%$}
& \multicolumn{5}{c}{$\alpha = 1\%$}
\\
\cline{2-16}
& \multicolumn{5}{c|}{}
& \multicolumn{5}{c|}{}
& \multicolumn{5}{c}{}
\\[-0.1in]
$n$ 
& $R$	&	$R^*_0$	&	$\overline{R}^*$	&	$\widehat{R}^*$	&	$\widetilde{R}^*$ 	
&	$R$	&	$R^*_0$	&	$\overline{R}^*$	&	$\widehat{R}^*$	&	$\widetilde{R}^*$ 	
&	$R$	& $R^*_0$	&	$\overline{R}^*$	& $\widehat{R}^*$	&	$\widetilde{R}^*$  \\
\hline            
40	&	12.0	&	10.8	&	10.2	&	9.0	&	9.5	&	6.9	&	5.8	&	5.5	&	4.3	&	4.9	&	1.9	&	1.4	&	1.3	&	0.8	&	1.0	\\
50	&	11.7	&	10.6	&	9.6	&	9.3	&	9.5	&	6.4	&	5.7	&	4.9	&	4.7	&	4.8	&	1.6	&	1.3	&	0.9	&	0.8	&	0.9	\\
60	&	11.6	&	10.9	&	10.1	&	9.5	&	9.9	&	6.3	&	5.7	&	4.9	&	4.7	&	5.0	&	1.3	&	1.1	&	0.9	&	0.7	&	0.8	\\
70	&	11.4	&	11.0	&	10.2	&	9.7	&	10.2	&	5.9	&	5.6	&	4.9	&	4.7	&	4.9	&	1.4	&	1.3	&	1.1	&	0.9	&	1.0	\\
100	&	10.7	&	10.3	&	9.7	&	9.7	&	9.7	&	5.9	&	5.6	&	5.0	&	5.1	&	5.2	&	1.2	&	1.1	&	0.9	&	0.9	&	1.0	\\
200	&	10.8	&	10.7	&	10.2	&	10.4	&	10.4	&	5.4	&	5.3	&	5.0	&	5.0	&	5.0	&	1.1	&	1.1	&	1.0	&	1.0	&	1.0	\\
\hline \hline
\end{tabular}
\end{table}
\begin{table}[!ht]
\begin{center}
\caption{Non-null rejection rates (\%), model 1, $n=40$, $\alpha=10\%$} \label{tab:sinalizada-linear-hetero-H0beta_poder}
\small
\begin{tabular}
{c| c c c c c c c c c c}
\hline \hline
$\epsilon$  &	0.1	&	0.5	&	1.0	&	2.0	&	4.0	&	6.0	\\
\hline
$R$	&	10.4	&	14.4	&	21.5	&	39.9	&	77.6	&	95.2	\\
$R^*_0$	&	10.3	&	14.7	&	21.9	&	40.5	&	78.2	&	95.8	\\
$\overline{R}^*$	&	10.3	&	14.6	&	21.7	&	40.0	&	78.3	&	95.8	\\
$\widehat{R}^*$	&	10.4	&	14.4	&	21.2	&	39.2	&	76.2	&	94.5	\\
$\widetilde{R}^*$ 	&	10.4	&	14.7	&	21.6	&	40.3	&	78.3	&	95.6	\\
\hline \hline
\end{tabular}
\end{center}
\end{table}

\begin{figure}[!ht]
  \centering
  \includegraphics[height=70mm,width=165mm]{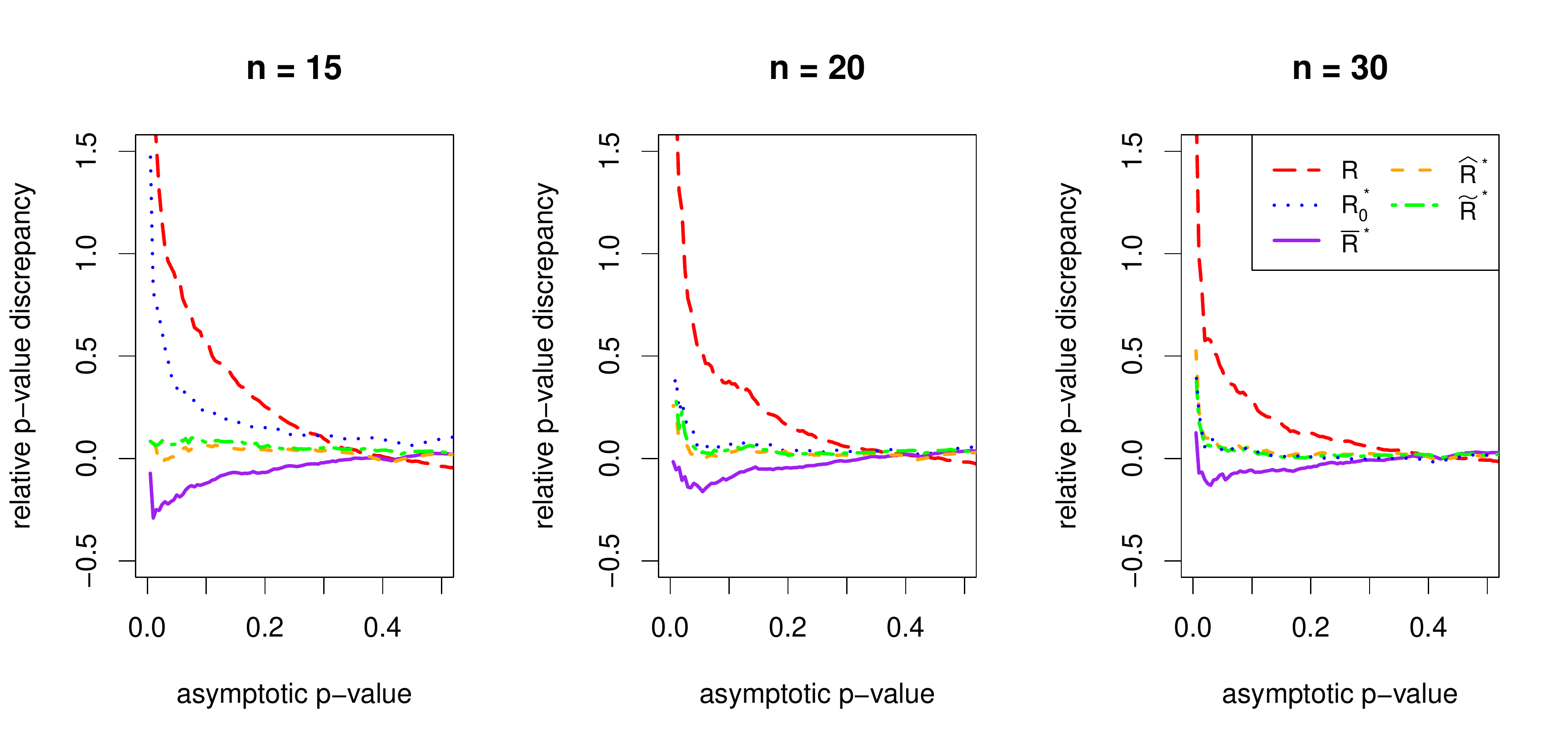}
  \caption{Relative p-value discrepancy plots, model 3.}
  \label{fig:sinalizada-naolinear-homo-valorp}
\end{figure}

We now consider the maximum extreme value regression model with constant dispersion and a nonlinear specification for the location parameter
given by 
$$
\mu_t = \beta_0 + \beta_1 x_{t1} +x_{t2}^ {\beta_2};
$$
`model 3'.
The null hypothesis is ${\mathcal H}_0: \beta_2 \leq 0$  and the alternative hypothesis is one-sided. 
We set $\sigma=1$, $\beta_0=1$, $\beta_1=1$ and $\beta_2=0$,
where the covariate values are randomly drawn from a ${\mathcal U}(0,1)$ distribution, and the sample sizes are 15, 20, 30 and 40.
The simulation results are presented in Tables~\ref{tab:sinalizada-naolinear-homo_tamanho} and \ref{tab:sinalizada-naolinear-homo_poder} and in Figure~\ref{fig:sinalizada-naolinear-homo-valorp}. 
As in the previous scenarios, the liberal behavior of the signed likelihood ratio test, the 
effectiveness of the adjustments and the similar power behavior of the different tests are evident. 
It should be noted, however, that the tests that use Severini's $\widehat{R}^*$ and Fraser, Reid and Wu's $\widetilde{R}^*$
statistics perform notably better than the other tests when the sample is small.

\begin{table}[htp]
\centering
\renewcommand{\arraystretch}{1.1}
\caption{Null rejection rates (\%); model 3}\label{tab:sinalizada-naolinear-homo_tamanho}
\footnotesize
\begin{tabular}
{c| c c c c c | c c c c c|  c c c c c}
\hline \hline
& \multicolumn{5}{c|}{$\alpha = 10\%$}
& \multicolumn{5}{c|}{$\alpha = 5\%$}
& \multicolumn{5}{c}{$\alpha = 1\%$}
\\
\cline{2-16}
& \multicolumn{5}{c|}{}
& \multicolumn{5}{c|}{}
& \multicolumn{5}{c}{}
\\[-0.1in]
$n$ 
& $R$	&	$R^*_0$	&	$\overline{R}^*$	&	$\widehat{R}^*$	&	$\widetilde{R}^*$ 	
&	$R$	&	$R^*_0$	&	$\overline{R}^*$	&	$\widehat{R}^*$	&	$\widetilde{R}^*$ 	
&	$R$	& $R^*_0$	&	$\overline{R}^*$	& $\widehat{R}^*$	&	$\widetilde{R}^*$  \\
\hline
15  &	14.2	&	12.0	&	8.7	&	10.6	&	10.8	&	8.1	&	6.6	&	3.9	&	5.1	&	5.5	&	2.2	&	1.7	&	0.7	&	1.1	&	1.1	\\
20	&	13.3	&	10.7	&	9.0	&	10.3	&	10.5	&	7.2	&	5.3	&	4.3	&	5.0	&	5.1	&	2.1	&	1.3	&	0.9	&	1.2	&	1.2	\\
30	&	12.1	&	10.4	&	9.4	&	10.5	&	10.4	&	6.8	&	5.3	&	4.6	&	5.3	&	5.2	&	1.8	&	1.2	&	1.0	&	1.2	&	1.2	\\
40	&	11.7	&	10.3	&	9.7	&	10.7	&	10.2	&	6.5	&	5.3	&	4.7	&	5.5	&	5.2	&	1.6	&	1.1	&	1.0	&	1.2	&	1.1	\\
100 & 10.5	&	10.0	&	9.7	&	10.0	&	10.0	&	5.9	&	5.4	&	5.1	&	5.4	&	5.4	&	1.2	&	1.0	&	0.9	&	1.1	&	1.0 \\
\hline \hline
\end{tabular}
\end{table}
\begin{table}[!ht]
\begin{center}
\caption{Non-null rejection rates (\%); model 3, $n=15$, $\alpha=10\%$} 
\label{tab:sinalizada-naolinear-homo_poder}
\small
\begin{tabular}
{c| c c c c c c c c c c}
\hline \hline
$\epsilon$	&	0.01	&	0.03	&	0.05	&	0.10	&	0.20	&	0.30	\\
\hline
$R$	&	12.0	&	16.2	&	21.1	&	37.5	&	76.5	&	96.3	\\
$R^*_0$	&	12.0	&	16.3	&	20.9	&	37.3	&	76.6	&	96.4	\\
$\overline{R}^*$	&	12.0	&	16.0	&	20.7	&	37.1	&	76.1	&	96.5	\\
$\widehat{R}^*$	&	11.8	&	16.3	&	21.0	&	36.8	&	75.5	&	95.6	\\
$\widetilde{R}^*$ 	&	12.0	&	16.2	&	20.9	&	37.3	&	76.5	&	96.3	\\
\hline \hline
\end{tabular}
\end{center}
\end{table}

\section{Applications}\label{applications}

Our first application involves with a data set (see Table \ref{tab:niwot}), consisting of maximum wind speed measured in January from 2001 to 2010 and the minimum temperature on the day in which the maximum wind speed was reached. 
The data were extracted from the alpine tundra climate station (3743 m) of Niwot Ridge/Green Lakes Valley, Colorado, USA.\footnote{CR23X data: {\tt http://niwot.colorado.edu/exec/.extracttoolA?d-1cr23x.ml}.}
\begin{table}[!ht]
\begin{center}
\caption{Maximum wind speed (m/s) and minimum temperature ($^{\rm o}$C) in January} 
\label{tab:niwot}
\begin{tabular}
{c| c c}
\hline\hline
year & temperature & wind speed\\
      \hline
2001 &    $-$7.40  &  33.42\\
2002 &   $-$11.95  &  44.04\\
2003 &   $-$17.99  &  42.92\\
2004 &   $-$25.63  &  42.51\\
2005 &   $-$16.61  &  45.75\\
2006 &   $-$10.93  &  47.78\\
2007 &    $-$9.21  &  43.34\\
2008 &   $-$26.13  &  48.69\\
2009 &   $-$20.27  &  43.20\\
2010 &   $-$19.00  &  43.00\\
\hline\hline
\end{tabular}
\end{center}
\end{table}

We assume that the maximum wind speed follows a maximum extreme value regression model (\ref{extreme value-max}) 
with constant dispersion and location sub-model given by
$$
\mu_t=\beta_0+\beta_1 x_{t},
$$
for $t=1,\ldots,10$, where the covariate $x$ is the minimum temperature. 

Table~\ref{tab:estimativas-niwot} presents the maximum likelihood estimates for the parameters and the corresponding standard errors.
The values of the test statistics for testing the null hypothesis ${\mathcal H}_0: \beta_1 \geq 0$ 
against ${\mathcal H}_1: \beta_1 \leq 0$ and the 
corresponding p-values are presented in Table~\ref{tab:niwot-sinalizada}. We note that the p-value of the unmodified
signed likelihood ratio test is 1.1\% while the p-values of the modified tests range from
2.8\% (Fraser, Reid and Wu's $\widetilde R^*$ and Barndorff-Nielsen's $R^*$) to 5.4\% (Skovgaard's $\overline R^*$).
The unmodified test displays the smallest p-value, in accordance with its liberal behavior observed in our simulation study.
The different adjustments weaken the evidence in favor of the alternative hypothesis.

\begin{table}[!ht]
\begin{center}
\caption{Estimates and standard errors, wind speed data set}
\label{tab:estimativas-niwot}
\begin{tabular}
{r r r r r r r r}
\hline \hline
Parameters    & $\beta_0$  & $\beta_1$       & $\sigma$ \\
\hline 
              & intercept  &  minimum        & \\
              &            &  temperature    & \\
\hline 
estimate &  34.3412  &  $-$0.4409 &   3.4211\\
s.e.     &   3.0910  &  0.1740  &   0.8435\\
\hline \hline
\end{tabular}
\end{center}
\end{table}

\begin{table}[!ht]
\begin{center}
\caption{Test statistics, , wind speed data set}\label{tab:niwot-sinalizada}
\small
\begin{tabular}
{r|r r r r r r}
\hline \hline 
            & $R$	&	$R^*_0$	&	$\overline{R}^*$	&	$\widehat{R}^*$	&	$\widetilde{R}^*$ &	$R^*$\\
\hline        
statistic	&	 $-$2.2912  &  $-$1.8989  &  $-$1.6085  &  $-$1.7592   &  $-$1.9043   &  $-$1.9043 \\
p-value   &    0.0110   &   0.0288    &  0.0539     & 0.0393     &  0.0284    &  0.0284  \\
\hline \hline 
\end{tabular}
\end{center}
\end{table} 

Our second application involves with a data set consisting of 34 men's decathlon results at the 1988 Olympic Games gathered by 
\citet[p. 304, data set 357]{HAND}.\footnote{The data set is
available at {\tt http://www.stat.ncsu.edu/working\_groups/sas/sicl/data/olympic.dat}.}
We assume that the score in high jump follows a maximum extreme value regression model (\ref{extreme value-max}) 
with constant dispersion and location sub-model 
$$
\mu_t=\beta_0+\beta_1x_{t1}+\beta_2x_{t2}+\beta_3x_{t3}+\beta_4x_{t4}+\beta_5x_{t5},
$$
for $t=1,\ldots,34$. The covariates are the scores in the following events: javelin throw ($x_1$), 
long jump ($x_2$), discus throw ($x_3$), shot put ($x_4$), and pole vault ($x_5$). 

The maximum likelihood estimates for the parameters and the corresponding standard errors are presented in Table~\ref{tab:estimativas-decatlo}, and the values of test statistics for testing the null hypotheses 
${\mathcal H}_0: \beta_1 \leq 0$ and ${\mathcal H}_0: \beta_4 \geq 0$ are presented in Table~\ref{tab:decatlo-sinalizada}.
For both cases, the p-values of the signed likelihood ratio test is much smaller than the p-values of the adjusted tests. 
It is interesting to observe that conclusions based on the unadjusted test may 
differ from those based on adjusted versions. 
For instance, at the 5\% nominal level, the signed likelihood ratio test ($R$) rejects ${\mathcal H}_0: \beta_1 \leq 0$, whereas all of the adjusted tests 
do not reject it. In addiction, at the 10\% nominal level, only the unadjusted test rejects the null hypothesis ${\mathcal H}_0: \beta_4 \geq 0$.
Based on our simulation results it is advisable to rely on the adjusted tests.

\begin{table}[!ht]
\begin{center}
\caption{Estimates and standard errors, decathlon data set}
\label{tab:estimativas-decatlo}
\begin{tabular}
{r r r r r r r r}
\hline \hline \\
Parameters  & $\beta_0$  & $\beta_1$ & $\beta_2$  & $\beta_3$  & $\beta_4$ & $\beta_5$  & $\sigma$ \\
\hline 
              & intercept &  javelin & long & discuss &  shot & pole  \\
              &           &    throw & jump &   throw &   put & vault \\
\hline 
estimate &  $-$0.1971  &  0.2125  &  0.0522  &  0.0371  & $-$0.1844  &  0.1522  &  0.0383 \\
s.e.     &      0.3873 &  0.0889  &  0.1674  &  0.1244  &    0.1291  &  0.0951  &  0.0051 \\
\hline \hline
\end{tabular}
\end{center}
\end{table}

\begin{table}[!ht]
\begin{center}
\caption{Test statistics, decathlon data set}\label{tab:decatlo-sinalizada}
\small
\begin{tabular}
{r|r r r r r r}
\hline \hline 
            & $R$	&	$R^*_0$	&	$\overline{R}^*$	&	$\widehat{R}^*$	&	$\widetilde{R}^*$ &	$R^*$\\
\hline        
${\mathcal H}_0: \beta_1 \leq 0$ versus ${\mathcal H}_1: \beta_1 > 0$ & & & & & &\\
 statistic	&	2.0102	&	1.5841	&	1.5658	&	1.5920	&	1.5945	&	1.5945\\
 p-value	    &	0.0222	&	0.0566	&	0.0587	&	0.0557	&	0.0554	&	0.0554\\
\hline
${\mathcal H}_0: \beta_4 \geq 0$ versus ${\mathcal H}_1: \beta_4 < 0$  & & & & & &\\
 statistic	&	$-$1.4815	&	$-$1.1743	&	$-$1.1277	&	$-$0.5062	&	$-$1.1789	&	$-$1.1789\\
 p-value	&	0.0692	&	0.1201	&	0.1297	&	0.3064	&	0.1192	&	0.1192\\
\hline \hline 
\end{tabular}
\end{center}
\end{table} 

\section{Conclusion}\label{conclusion}

In this paper we derived Barndorff-Nielsen's modified signed likelihood ratio statistic 
for the homoskedastic linear extreme value regression model, and different approximations 
to this statistic in the general case were made, i.e., possibly nonlinear and/or heteroskedastic extreme value 
regression models. Unlike the other adjusted statistics, DiCiccio \& Martin's statistic requires an 
orthogonal parameterization, which makes its deduction more involved.
Our simulation results revealed that the signed likelihood ratio test is liberal in small and moderate-sized samples, 
i.e., it leads to a type I error probability greater than the chosen nominal level. Additionally, it is evident from 
the simulations that all of the adjustments considered in this paper are effective in shrinking the size distortion 
of the original test. For the linear homoskedastic model, Barndorff-Nielsen's, DiCiccio's and Fraser, Reid \& Wu's tests
behaved equally well and clearly better than Severini's and Skovgaard's tests. For the nonlinear homoskedastic model,
Severini's and Fraser, Reid \& Wu's tests presented better performance than the others. 
In terms of power, all of the tests behaved similarly.
Overall, Fraser, Reid \& Wu's test is the best performing test. 
 
Our applications show that the conclusion reached by the signed likelihood ratio test may conflict with the conclusion achieved by the adjusted tests when the sample is not large (in our application, $n=10$ and 34). As indicated by our simulations, practitioners should rely on the adjusted tests.

\section*{Acknowledgments}
We gratefully acknowledge the financial support from the Brazilian agencies CNPq, CAPES and FAPESP.

\small
\appendix
\section*{Appendix}
\label{ape:appendix}

Let $v$ be a vector, and $v_{(r)}$ be the vector $v$ without its $r$th component.
Analogously, let $\mathcal V$ be a matrix and $\mathcal V_{r}$ be the matrix $\mathcal V$
without the $r$th column. 
Assume that $\nu=\beta_r$ is the parameter of interest. 
The Fisher information matrix (see \citealp{FERRARI2012}) can be written as
$$
I(\theta)=I(\beta_r,\theta_{(r)})=I(\beta_r,\beta_{(r)},\gamma)=
\begin{pmatrix}
I_{\beta_r\beta_r} & I_{\beta_r\beta_{(r)}} & I_{\beta_r\gamma} \\
I_{\beta_{(r)}\beta_r} & I_{\beta_{(r)}\beta_{(r)}} & I_{\beta_{(r)}\gamma} \\
I_{\gamma\beta_r} & I_{\gamma\beta_{(r)}} & I_{\gamma\gamma} \\
\end{pmatrix} ,
$$
where
$$
\begin{array}{l l}
I_{\beta_r\beta_r}={\mathbf x}_{\cdot r}^{\!\top}\Phi^{-1}T^2\Phi^{-1}{\mathbf x}_{\cdot r},
&
I_{\beta_r\beta_{(r)}}=I_{\beta_{(r)}\beta_r}^{\!\top}={\mathbf x}_{\cdot r}^{\!\top}\Phi^{-1}T^2\Phi^{-1}X_{(r)},
\\
I_{\beta_r\gamma}=I_{\gamma\beta_r}^{\!\top}=({\mathcal E}-1){\mathbf x}_r^{\!\top}\Phi^{-1}TH\Phi^{-1}Z,
&
I_{\beta_{(r)}\beta_{(r)}}=X_{(r)}^{\!\top}\Phi^{-1}T^2\Phi^{-1}X_{(r)},
\\
I_{\beta_{(r)}\gamma}=I_{\gamma\beta_{(r)}}^{\!\top}=({\mathcal E}-1)X_{(r)}^{\!\top}\Phi^{-1}TH\Phi^{-1}Z,
&
I_{\gamma \gamma }=\bigl(1 + \Gamma ^{(2)}(2)\bigr)Z^{\!\top}\Phi^{-1}H^2\Phi^{-1}Z,
\end{array}
$$
${\mathbf x}_{\cdot r}=(x_{1r},\ldots,x_{nr})^{\!\top}$ is the $r$th column of $X$ 
and $\Phi$, $T$, $H$, $X$ and $Z$ as defined above.

Consider a parameterization $\vartheta=(\beta_r,\vartheta_{(r)})$, $\vartheta_{(r)}=(\kappa,\tau)$, where 
$\vartheta$, $\vartheta_{(r)}$, $\kappa$ and $\tau$ have dimensions $(k+m)$, $(k+m-1)$, $(k-1)$ and $m$, respectively, in such a way 
that $\beta_r$ is orthogonal to $\vartheta_{(r)}$. 
Let
$$
I_{\beta_r\theta_{(r)}}=
\begin{pmatrix}
I_{\beta_r\beta_{(r)}} & I_{\beta_r\gamma}
\end{pmatrix},
\qquad
I_{\theta_{(r)}\theta_{(r)}}=
\begin{pmatrix}
I_{\beta_{(r)}\beta_{(r)}} & I_{\beta_{(r)}\gamma}\\
I_{\gamma\beta_{(r)}} & I_{\gamma\gamma}
\end{pmatrix},
\qquad
I_{\theta_{(r)}\theta_{(r)}}^{-1}=
\begin{pmatrix}
I^{\beta_{(r)}\beta_{(r)}} & I^{\beta_{(r)}\gamma}\\
I^{\gamma\beta_{(r)}} & I^{\gamma\gamma}
\end{pmatrix}.
$$
Define
$
A \equiv I_{\beta_r\beta_{(r)}}I^{\beta_{(r)}\beta_{(r)}}+I_{\beta_r\gamma}I^{\gamma\beta_{(r)}},
$
$
B \equiv I_{\beta_r\beta_{(r)}}I^{\beta_{(r)}\gamma}+I_{\beta_r\gamma}I^{\gamma\gamma},
$
$\beta_{(r)}=\kappa-\beta_rA^{\!\top}$ and $\gamma=\tau-\beta_rB^{\!\top}$. 
From \citet{COX1987}, it can be shown that $\beta_r$ is orthogonal to $(\kappa, \tau)$.

In order to obtain $R_0^*$, one should first re-write (\ref{extreme value-max})-(\ref{linkphi}) in the orthogonal parameterization  $\vartheta=(\beta_r,\kappa, \tau)$, and define the systematic components $\eta^*(\vartheta)$ and $\delta^*(\vartheta)$. The derivative of the log-likelihood function with respect to the parameter of interest $\beta_r$, $\ell_{\beta_r}^*(\vartheta)$, the observed information matrix, $J^*(\vartheta)$ and the element of the expected information matrix that corresponds to $\beta_r$, $I_{\beta_r\beta_r}^*(\vartheta)$, should then be computed in this parameterization as in \citet[p. 584]{FERRARI2012}, and inserted into (\ref{U0}).

Consider, for instance, the maximum extreme value model regression model (\ref{extreme value-max}) with location and dispersion sub-models
$
g(\mu_t) = \eta_t = \eta(x_{t},\beta)=x_{t}^\top\beta \quad {\rm and}
$
and
$
h(\sigma_t) =  \delta_t = \delta(z_{t},\gamma)=z_{t}^\top\gamma.
$
We have
\begin{equation}\label{etareparameterization}
g(\mu_t)=
x_{tr}\beta_r+x_{t(r)}^\top\beta_{(r)}=
x_{tr}\beta_r+x_{t(r)}^\top(\kappa-\beta_r A^{\!\top})=
(x_{tr}-x_{t(r)}^\top A^{\!\top})\beta_r+x_{t(r)}^\top\kappa \equiv \eta^*_t(\vartheta)
\end{equation} 
and
\begin{equation}\label{deltareparameterization}
h(\sigma_t)=
z_{t}^\top \tau-\beta_r z_{t}^\top B^{\!\top}\equiv\delta^*_t(\vartheta).
\end{equation}
For the reparameterized model (\ref{extreme value-max}) with (\ref{etareparameterization})-(\ref{deltareparameterization}),
we have
$$
\ell^*_{\beta_r}(\vartheta)=\frac{\partial \ell^*(\vartheta)}{\partial\beta_r}=
\iota^{\!\top}\Bigl ( \Phi ^{-1}T\left( {\mathcal I} - \breve {\mathcal Z} \right) {\mathbf v}_1
+
\Phi^{-1}H\Bigl(-{\mathcal I}+{\mathcal Z}-{\mathcal Z}\breve {\mathcal Z} \Bigr){\mathbf v}_2   \Bigr),
$$
$$
\begin{array}{l l}
J^*_{\beta_r\beta_r}= 
{\mathbf v}_1^{\!\top}   V_{\beta\beta} {\mathbf v}_1 
+2{\mathbf v}_1^{\!\top}V_{\beta\gamma} \; {\mathbf v}_2
+{\mathbf v}_2^{\!\top}V_{\gamma\gamma} \; {\mathbf v}_2,
&
J^*_{\kappa \kappa }=X_{(r)}^{\!\top} \; V_{\beta\beta} \; X_{(r)}, \\
J^*_{\beta_r\kappa}=J_{\kappa\beta_r}^{\!\top}= 
{\mathbf v}_1^{\!\top}   V_{\beta\beta} \; X_{(r)}+{\mathbf v}_2^{\!\top} \; V_{\beta\gamma} \; X_{(r)},
&
J^*_{\kappa\tau}=J_{\tau\kappa}^{\!\top}=X_{(r)}^{\!\top}V_{\beta\gamma} \; Z, \\
J^*_{\beta_r\tau}=J_{\tau\beta_r}^{\!\top}= 
{\mathbf v}_1^{\!\top}V_{\beta\gamma} \; Z 
+{\mathbf v}_2^{\!\top}V_{\gamma\gamma} \; Z,
&
J^*_{\tau \tau }= Z^{\!\top}V_{\gamma\gamma} \; Z
\end{array}
$$
and
$$
I^*_{\beta_r\beta_r}= {\mathbf v}_1^{\!\top}   \Phi ^{-1}T^2\Phi ^{-1} {\mathbf v}_1 
+2({\mathcal E}-1){\mathbf v}_1^{\!\top}   \Phi ^{-1}T H \Phi ^{-1} {\mathbf v}_2 
+\bigl(1+\Gamma^{(2)}(2)\bigr){\mathbf v}_2^{\!\top}\Phi^{-1}H^2\Phi^{-1}{\mathbf v}_2,
$$
where
${\mathbf v}_1=\left( {\mathbf x}_{\cdot r}-X_{(r)}A^{\!\top}  \right)$,
${\mathbf v}_2=-(ZB^{\!\top})$,
$V_{\beta\beta}=\Phi ^{-1}T\Bigl(\breve {\mathcal Z}\Phi ^{-1} + ({\mathcal I}-\breve {\mathcal Z})ST\Bigr)T$,
$V_{\beta\gamma}=\Phi^{-1}T({\mathcal I}-\breve {\mathcal Z}+{\mathcal Z}\breve {\mathcal Z})H\Phi^{-1}$ 
and
$V_{\gamma\gamma}=\Phi^{-1}H\Bigl((-{\mathcal I}+2{\mathcal Z}-2{\mathcal Z}\breve {\mathcal Z}+{\mathcal Z}^2\breve {\mathcal Z})\Phi^{-1}+(-{\mathcal I}+{\mathcal Z}-{\mathcal Z}\breve {\mathcal Z})QH\Bigr)H$. 
DiCiccio \& Martin's adjusted signed likelihood ratio statistic can now be obtained by replacing the formulas above 
in (\ref{U0}).

\end{document}